\documentclass[letterpaper]{amsart}
\usepackage{hyperref}
\usepackage{amsrefs}
\usepackage{amsthm}
\usepackage{amssymb}
\usepackage{amsmath}
\usepackage{xcolor}
\usepackage{easyReview}

\numberwithin{equation}{section}

\newtheorem{thm}{Theorem}

\newtheorem{prop}[thm]{Proposition}

\theoremstyle{definition}

\theoremstyle{remark}
\newtheorem{rmk}[thm]{Remark}

\begin{document}

\title[Revisiting Jacobi--Trudi identities via the BGG category $\mathcal{O}$]{Revisiting Jacobi--Trudi identities via the BGG category $\mathcal{O}$}

\author{Tao Gui}
\address{Beijing International Center for Mathematical Research, Peking University, Beijing, 100871, P.R. China}
\email{guitao18@mails.ucas.ac.cn}

\author{Arthur L. B. Yang}
\address{Center for Combinatorics, LPMC, Nankai University, Tianjin 300071, P. R. China}
\email{yang@nankai.edu.cn}

\begin{abstract}
By interpreting Kostka numbers as tensor product multiplicities in the BGG category $\mathcal{O}$ for the special linear Lie algebras, we provide a new proof of the classical Jacobi--Trudi identities for skew Schur polynomials, derived from the celebrated Weyl character formula. We re-establish the Schur positivity of certain truncations in the Jacobi--Trudi expansion of skew Schur polynomials and obtain Schur positivity results for similar truncations in the Jacobi--Trudi-type expansion of the product of two Schur polynomials. Furthermore, we interpret the coefficients in the Schur polynomial expansions of these Jacobi--Trudi truncations as tensor product multiplicities in the BGG category $\mathcal{O}$. 
\end{abstract}


\maketitle

\section{Introduction}

Many open problems in symmetric function theory are related to the \emph{Schur positivity} questions, that is, showing that certain symmetric functions can be expressed as non-negative linear combinations of \emph{Schur functions}; see, for example, \cites{haiman1993jams,stanley1998graph,haiman2001jams, lam2007schur,chari2014posets}.
For our purpose here, we only need to consider \emph{Schur polynomials}, which are obtained from Schur functions by reducing the number of variables.
The significance of Schur positivity stems from the fact that Schur polynomials not only have interesting combinatorial interpretations but also admit natural algebraic or geometric meanings. In combinatorics, Schur polynomials can be defined as generating functions of weighted semi-standard Young tableaux. In representation theory, they appear as characters of finite-dimensional irreducible polynomial representations (known as \emph{Schur modules}) of the general linear groups $GL_{n}(\mathbb{C})$, as well as Frobenius images of irreducible representations of the symmetric group. The relationship between Schur polynomials and representation theory was discovered by Schur, hence these important polynomials are named after him. In geometry, Schur polynomials appear as the cohomology classes of the Schubert cycles of the corresponding Grassmannian. For more information on Schur polynomials, we refer the reader to \cites{fulton1997young,Stanley99Enumerative,macdonald1998symmetric}. Due to these connections, whether a particular symmetric polynomial can be expanded positively in terms of Schur polynomials is of particular interest in symmetric function theory, representation theory, and classical Schubert calculus. 

The main objective of this paper is to study the Schur positivity of certain truncations of the Jacobi--Trudi expansion of skew Schur polynomials and the Jacobi--Trudi-type expansion of the product of two Schur polynomials. For ordinary partitions, this positivity can be implied by the existence of a conjectural resolution (due to Lascoux \cite{lascoux1977polynomes}) of irreducible polynomial representations of $GL_{n}(\mathbb{C})$, which was later constructed by Akin\cite{akin1988complexes} and Zelevinskii\cite{zelevinskii1987resolvents} independently. Zelevinskii's approach could also handle the case of skew partitions. However, it is an open problem to give a combinatorial interpretation of these non-negative coefficients in these Schur polynomial expansion\cite{Stanley99Enumerative}*{p.461}. 
By the work of Akin and Zelevinskii, one can interpret these coefficients as multiplicities of irreducible representations of some polynomial representation of $GL_{n}(\mathbb{C})$. Using the classical BGG resolution and techniques from the BGG category $\mathcal{O}$ of the special linear Lie algebras $\mathfrak{s l}_{n}(\mathbb{C})$, we can prove the Schur positivity for the case of skew partitions and for the case of products. We interpret these coefficients as tensor product multiplicities in the BGG category $\mathcal{O}$ of $\mathfrak{s l}_{n}(\mathbb{C})$, thus providing a new algebraic interpretation.  

In order to state our results explicitly, let us now recall some related definitions and results. We will adopt the terminology and notation of symmetric functions in \cite{Stanley99Enumerative}. An integer partition $\lambda$ of a positive integer $d$, denoted by $\lambda\vdash d$,  is a weakly decreasing sequence $(\lambda_1,\lambda_2,\ldots,\lambda_n)$ of non-negative integers such that 
$|\lambda|:=\sum_{i=1}^n\lambda_i=d$. For each partition $\lambda$, let $m_{\lambda},\, h_{\lambda},\, s_{\lambda}$ denote the corresponding monomial symmetric polynomial, the complete symmetric polynomial, and the Schur polynomial, respectively.
It is known that each of $\{m_{\lambda}\,|\,\lambda\vdash d\}$, $\{h_{\lambda}\,|\,\lambda\vdash d\}$, and $\{s_{\lambda}\,|\,\lambda\vdash d\}$ is a basis of the space of homogeneous symmetric polynomials of degree $d$. Here we adopt the definition of the Schur polynomial $s_{\lambda}$ given by 
\begin{equation} \label{def}
s_{\lambda}:=\sum_{\mu} K_{\lambda, \mu} m_{\mu},
\end{equation}
where each $K_{\lambda, \mu}$, called the \emph{Kostka number}, is equal to the number of semistandard Young tableaux of shape $\lambda$ and type $\mu$.

If $\mu,\nu$ are two partitions of length at most $n$ such that $\mu_i\geq \nu_i$ for each $i$, then one can define a skew partition $\mu/\nu$.
Let $s_{\mu/\nu}$ denote the corresponding skew Schur polynomial indexed by $\mu/\nu$, which can be defined by the property
\begin{equation} \label{skew}
\left\langle s_{\mu / \nu}, s_{\lambda}\right\rangle=\left\langle s_{\mu}, s_{\lambda} s_{\nu}\right\rangle,
\end{equation}
where the scalar product is defined by claiming that Schur polynomials form an orthonormal basis. The classical \emph{Jacobi--Trudi identity} states that
\begin{equation} \label{GJK}
s_{\mu / \nu}=\operatorname{det}\left(h_{\mu_{i}-\nu_{j}-i+j}\right)_{i, j=1}^{n},
\end{equation}
where we set $h_{0}=1$ and $h_{k}=0$ for $k<0$. There is also a Jacobi--Trudi-type identity for the product of two Schur polynomials
:
\begin{equation} \label{eq-JT type}
   s_\mu s_\nu=\operatorname{det}\left(h_{\mu_i+\nu_{n+1-j}-i+j}\right)_{i, j=1}^n, 
\end{equation}
see \cite{macdonald1998symmetric}*{Examples 3.8, pp. 46-47}.

There are several proofs of the Jacobi--Trudi identity. The classical combinatorial proof is by relating the semi-standard Young tableaux with non-intersection lattice paths and then using the Lindstr\"om--Gessel--Viennot lemma; see \cite{sagan2001symmetric}*{Section 4.5}, or the first proof of \cite{Stanley99Enumerative}*{Theorem 7.16.1}. A purely algebraic proof can also be given; see \cite{macdonald1998symmetric}*{Section 1.3}, or the second proof of \cite{Stanley99Enumerative}*{Theorem 7.16.1}. Since Schur polynomials also appear as characters of Schur modules of the general linear groups, the alternating sum format of the spanning of the right hand side determinant in \eqref{GJK} resembles an \emph{Euler--Poincar\'e characteristic}, so it is natural to ask whether there is a resolution of Schur modules to realize the identity \eqref{GJK} at least for the case of $\nu=\emptyset$. The existence of such a resolution was conjectured by Lascoux \cite{lascoux1977polynomes} and was later constructed by Akin in \cites{akin1988complexes,akin1989resolutions,akin1992complexes} and independently by Zelevinskii in \cite{zelevinskii1987resolvents}. 

The existence of such a resolution has the following implication, see \cite[Exercise 7.38]{Stanley99Enumerative}.

\begin{thm}[Akin--Zelevinskii theorem] 
For $0 \leq k \leq \binom{n}{2}$, the symmetric polynomial
\begin{equation} \label{truncated}
g_{\mu / \nu}^{k}=(-1)^{k} \sum_{w \in S_{n}, \ell(w) \geq k} (-1)^{\ell(w)} h_{\mu+\delta-w(\nu+\delta)},
\end{equation}
a truncation of the Jacobi--Trudi expansion of $s_{\mu/\nu}$,  
is Schur positive, where $\mu$ and $\nu$ are partitions of length at most $n$, $\ell(w)$ is the number of inversions of the permutation $w \in S_{n}$, and $\delta=(n-1, n-2, \ldots, 1,0)$.
\end{thm}

In this paper, We provides a new proof of the above theorem via the BGG category $\mathcal{O}$ of $\mathfrak{s l}_{n}(\mathbb{C})$. Using the duality of finite dimensional simple modules, we prove the following result.

\begin{thm}
With the same notations as in the above theorem, for $0 \leq k \leq \binom{n}{2}$, the symmetric polynomial
\begin{equation} 
t_{\mu,\nu}^{k}=(-1)^{k} \sum_{w \in S_{n}, \ell(w) \geq k} (-1)^{\ell(w)} h_{\mu+\delta+w(w_0 \nu-\delta)},
\end{equation}
a truncation of the following Jacobi--Trudi-type expansion of $s_{\mu}s_{\nu}$,
is Schur positive, where $w_0$ is the longest element of $S_{n}$. 
\end{thm}

Actually, we give a following new interpretation of the coefficients in the Schur polynomial expansion of $g_{\mu / \nu}^{k}$ and $t_{\mu, \nu}^{k}$ via the BGG category $\mathcal{O}$ of $\mathfrak{s l}_{n}(\mathbb{C})$, which in turn leaves the above two theorems as immediate corollaries.

\begin{thm} \label{thm1}
For $0 \leq k \leq \binom{n}{2}$, we have
\begin{equation} \label{new}
g_{\mu / \nu}^{k}=\sum_{\lambda}
[L(\lambda) \otimes \operatorname{V}\left(\nu, k\right): L(\mu)] s_{\lambda},
\end{equation}
and
\begin{equation}
\label{eq-tmunu}
t_{\mu,\nu}^{k}=\sum_{\lambda}
[L(\lambda) \otimes \operatorname{V}\left(-w_0\nu, k\right): L(\mu)] s_{\lambda},
\end{equation}
where $L(\lambda)$ (and $L(\mu)$, respectively) denotes the simple module of highest weight $\lambda$ (and $\mu$, respectively) of the special linear Lie algebra $\mathfrak{s l}_{n}(\mathbb{C})$, $\operatorname{V}\left(\nu, k\right)$ is the image of the boundary map $d_{k}$ in the BGG resolution of $L(\nu)$:
\begin{equation}
0 \rightarrow C_{\binom{n}{2}}^{\nu} \stackrel{d_{\binom{n}{2}}}{\longrightarrow} C_{\binom{n}{2}-1}^{\nu} \rightarrow \cdots \rightarrow C_{1}^{\nu} \stackrel{d_1}{\longrightarrow} C_{0}^{\nu} \stackrel{d_0}{\longrightarrow} L(\nu) \rightarrow 0,
\end{equation}
$\operatorname{V}\left(-w_0\nu, k\right)$ is similarly the image of a boundary map in the BGG resolution of $L(-w_0\nu)$ and the bracket denotes the composition factor multiplicity in the BGG category $\mathcal{O}$.
\end{thm}

Our approach is different in the sense that we transform the Schur positivity problem into the positivity problem of a virtual character in the BGG category $\mathcal{O}$ and then use classical BGG resolution to prove the positivity. It also allows us to give a new proof of the Jacobi--Trudi identity \eqref{GJK}. Actually, we can derive it from the Weyl character formula, as illustrated in subsequent sections. 

We observed that truncations of the Jacobi--Trudi expansion of a skew Schur polynomial seem to share similar convexity of Schur polynomials, see \cite{huh2022logarithmic}. Specifically, with Candice X. T. Zhang, we have checked that for partitions $\mu$ and $\nu$ with $|\mu|+|\nu| \leq 9$, all \emph{normalized} truncations of the Jacobi--Trudi expansion of a skew Schur polynomial $g_{\mu / \nu}^{k}$ are \emph{Lorentzian}. It would be interesting to see whether our approach can be used to understand the Lorentzian property for the skew Schur polynomials, see \cite{huh2022logarithmic}*{Conjecture 19}.

The remaining part of this paper is organized as follows. In Section \ref{sect-prel}, we provide preliminary results on the BGG category $\mathcal{O}$ and the classical BGG resolution. In Section \ref{sect-main}, we first give a new interpretation of  Kostka numbers in the BGG category $\mathcal{O}$, and then leave the world of symmetric functions and use Verma module characters to prove the Jacobi--Trudi identity and Theorem \ref{thm1}. 

\setcounter{tocdepth}{2}

\section{Preliminaries}\label{sect-prel}

In this section, we recall the relevant representation theory of complex semi-simple Lie algebras. For a more detailed background,  we refer the reader to the comprehensive reference \cite{humphreys2008representations}. 

\subsection{The BGG category $\mathcal{O}$}

Let $\mathfrak{g}$ be a finite-dimensional complex semi-simple Lie algebra (for our application, $\mathfrak{g}=\mathfrak{sl}_{n}(\mathbb{C})$). We fix a Cartan decomposition $\mathfrak{g}=\mathfrak{n}_{-} \oplus \mathfrak{h} \oplus \mathfrak{n}_{+}$. (In the case of $\mathfrak{g}=\mathfrak{s l}_{n}(\mathbb{C})$, we could choose $\mathfrak{n}_{+}$ and $\mathfrak{n}_{-}$ to be the strictly upper-triangular and strictly lower-triangular matrices, respectively, and $\mathfrak{h}$ to be the diagonal matrices.) We consider the \emph{Bernstein--Gelfand--Gelfand category $\mathcal{O}$}, which is a full subcategory of $\operatorname{Mod} U(\mathfrak{g})$ (the category of left-modules of the enveloping algebra of $\mathfrak{g}$). Each object of $\mathcal{O}$ is a module $M$ satisfying the following three conditions:
\begin{itemize}
\item[($\mathcal{O} 1$)] $M$ is a finitely generated $U(\mathfrak{g})$-module.

\item[($\mathcal{O} 2$)] $M$ is $\mathfrak{h}$-semi-simple, that is, $M$ is a weight module: $M=\bigoplus_{\lambda \in \mathfrak{h}^{*}} M_{\lambda}$, where $M_\lambda=\{m \in M: h \cdot m=\lambda(h) m$ for all $h \in \mathfrak{h}\}$.

\item[($\mathcal{O} 3$)] $M$ is locally $\mathfrak{n}_{+}$-finite, that is, for each $v \in M$, the subspace $U(\mathfrak{n}_{+}) \cdot v$ of $M$ is finite-dimensional.
\end{itemize}

For example, all finite-dimensional modules lie in $\mathcal{O}$. It follows from the axioms that all weight spaces of $M \in \mathcal{O}$ are finite-dimensional. Thus one could define the \emph{(formal) character} of $M \in \mathcal{O}$ by
\begin{align*}
\operatorname{ch}(M)=\sum_{\lambda} \operatorname{dim}\left(M_{\lambda}\right) x^{\lambda},
\end{align*}
where the symbols $x^{\lambda}$ denote formal exponentials, multiplying by the rule $x^{\lambda}  x^{\mu}=x^{\lambda+\mu}$.
Furthermore, if $M \in \mathcal{O}$ and $L$ is finite-dimensional, then $L \otimes M \in \mathcal{O}$ and $\operatorname{ch}(L \otimes M)=\operatorname{ch} L  \operatorname{ch} M$.

For each $\lambda \in \mathfrak{h}^{*}$, there is a one-dimensional representation $\mathbb{C}_{\lambda}$ of $\mathfrak{b}:=\mathfrak{h} \oplus \mathfrak{n}_{+}$, where $\mathfrak{h}$ acts via $\lambda$ and $\mathfrak{n}_{+}$ acts by zero. Now consider the \emph{Verma (or standard) module} $\Delta(\lambda)$, defined as  $\Delta(\lambda):=U(\mathfrak{g}) \otimes_{U(\mathfrak{b})} \mathbb{C}_{\lambda}$,
where $U(\mathfrak{g})$ and $U(\mathfrak{b})$ are the universal enveloping algebras of $\mathfrak{g}$ and $\mathfrak{b}$, respectively.
It can be proved that simple objects $L(\lambda)$ (\emph{simple highest weight module}) in $\mathcal{O}$ are parametrized by $\lambda \in \mathfrak{h}^{*}$ and can be uniformly constructed as the unique simple quotient of the corresponding $\Delta(\lambda)$. 

Since each $M \in \mathcal{O}$ is both artinian and noetherian, it follows that $M$ possesses a composition series with simple quotients isomorphic to various $L(\lambda)$. The \emph{multiplicity} of $L(\lambda)$ in $M$ is independent of the choice of composition series and is denoted by $[M: L(\lambda)]$. It is known that category $\mathcal{O}$ is an abelian category and if $0 \rightarrow M^{\prime} \rightarrow M \rightarrow M^{\prime \prime} \rightarrow 0$ is a short exact sequence in $\mathcal{O}$, then $\operatorname{ch} M=\operatorname{ch} M^{\prime}+\operatorname{ch} M^{\prime \prime}$. Thus any $\operatorname{ch} M$ is determined uniquely by the formal characters of the composition factors of $M$, taken with multiplicity.

\subsection{Weyl character formula and BGG resolution}

One of the basic questions in the representation theory of the complex semi-simple Lie algebra $\mathfrak{g}$ is to compute the character of $L(\lambda)$.
Those $L(\lambda)$ with $\lambda$ dominant and integral are precisely all the finite-dimensional simple modules of $\mathfrak{g}$, and the characters are given by the celebrated Weyl character formula
\footnote{It is easy to see that Cauchy's bialternant formula for the Schur polynomial $s_{\lambda}\left(x_{1}, x_{2}, \ldots, x_{n}\right)=\operatorname{det}\left(x_{j}^{\lambda_{i}+n-i}\right) / \operatorname{det}\left(x_{j}^{n-i}\right)$, where $\lambda=\left(\lambda_{1}, \lambda_{2}, \ldots, \lambda_{n}\right)$ is a partition, is a special case of the Weyl character formula for $\mathfrak{g}=\mathfrak{s l}_{n}(\mathbb{C})$ and $W=S_{n}$.}

\begin{equation} \label{Weyl}
\operatorname{ch}L(\lambda)=\frac{\sum_{w \in W} (-1)^{\ell(w)} x^{w(\lambda+\rho)}}{\sum_{w \in W} (-1)^{\ell(w)} x^{w(\rho)}},
\end{equation}
where $W$ is the Weyl group, and $\rho$ is the Weyl vector, defined as half the sum of positive roots. For $\mathfrak{g}=\mathfrak{s l}_{n}(\mathbb{C})$, the Weyl group $W$ is the permutation group $S_{n}$, and $\rho$ corresponds to the partition $\delta=(n-1, n-2, \ldots, 1,0)$. It can be implied from the Weyl character formula that the dual space $L(\lambda)^*$, with the standard action $(x \cdot f)(v)=-f(x \cdot v)$ for $x \in \mathfrak{g}, v \in L(\lambda), f \in L(\lambda)^*$, is isomorphic to $L\left(-w_0 \lambda\right)$ (where $w_{0} \in W$ is the longest element).

In terms of characters of Verma modules, the Weyl character formula is equivalent to
\begin{equation} \label{WK}
\operatorname{ch} L(\lambda)=\sum_{w \in W} (-1)^{\ell(w)} \operatorname{ch} \Delta(w \cdot \lambda),
\end{equation}
where the \emph{dot action} is defined as $w \cdot \lambda =w(\lambda+\rho)-\rho$, in other words, this shifts the standard action of $W$ to have center $-\rho$, and the characters of the Verma modules are given by 
\begin{align}\label{chac-verma}
\operatorname{ch}\left(\Delta(\lambda)\right)=x^{\lambda} / \prod\left(1-x^{-\alpha}\right),
\end{align}
where $\alpha$ ranges over all positive roots of the complex semi-simple Lie algebra $\mathfrak{g}$. 

In a classical paper \cite{bernstein1975differential},  Bernstein, Gelfand and Gelfand studied the Weyl character formula and constructed a resolution of $L(\lambda)$ for $\lambda$ dominant and integral to realize the alternating sum in the Weyl character formula (\ref{WK}) as an Euler--Poincar\'e characteristic of a complex, which is now known as the \emph{BGG resolution}. Specifically, they proved the following theorem.
\begin{thm} \cite{bernstein1975differential}*{Theorem 10.1} \label{thm4}
Let $\lambda$ be a dominant and integral weight and $L(\lambda)$ be the irreducible finite-dimensional $\mathfrak{g}$-module with highest weight $\lambda$. Then there exists an exact sequence of $\mathfrak{g}$-modules 
\begin{equation} \label{BGG}
0 \rightarrow C_{\binom{n}{2}}^{\lambda} \stackrel{d_{\binom{n}{2}}}{\longrightarrow} C_{\binom{n}{2}-1}^{\lambda} \rightarrow \cdots \rightarrow C_{1}^{\lambda} \stackrel{d_1}{\longrightarrow} C_{0}^{\lambda} \stackrel{d_0}{\longrightarrow} L(\lambda) \rightarrow 0,
\end{equation}
where $C_{k}^{\lambda}=\bigoplus_{w \in W, \ell(w)=k} \Delta(w\cdot\lambda)$.
\end{thm}

For general $\lambda \in \mathfrak{h}^{*}$, it is not hard \footnote{However, the determination of the coefficients $b(\lambda, w)$ is rather difficult, and they are given by the famous Kazhdan--Lusztig polynomials in virtue of the Kazhdan--Lusztig conjecture. We do not need the exact formula in this paper.} to show that
\begin{equation} \label{uni}
\operatorname{ch} L(\lambda)=\sum_{w\cdot\lambda \leq \lambda} b(\lambda, w) \operatorname{ch} \Delta(w\cdot\lambda) \text { with } b(\lambda, w) \in \mathbb{Z} \text { and } b(\lambda, \lambda)=1,
\end{equation}
where the partial ordering is defined by $\mu \leq \lambda$ if and only if $\lambda-\mu$ is a non-negative linear combination of simple roots. That is, $b(\lambda, w)$ is uni-triangular with respect to the partial ordering. It follows that formal characters $\operatorname{ch} L(\lambda)$ with $\lambda \in \mathfrak{h}^{*}$ are linearly independent and form a $\mathbb{Z}$-basis of the additive group generated by all $\operatorname{ch} \Delta(\lambda)$, while formal characters of Verma modules provide an alternative $\mathbb{Z}$-basis. The multiplicity $[M: L(\lambda)]$ of $L(\lambda)$ in $M$ is equal to the coefficient of $\operatorname{ch} L(\lambda)$ in the expansion of $\operatorname{ch} M$ in terms of irreducible characters.

\section{Proofs of the Main results}\label{sect-main}

The purpose of this section is two-fold: first, we give a new proof of the Jacobi--Trudi identity, and secondly, we use BGG resolution to prove Theorem \ref{thm1}. 

Our approach given here relies on an interpretation of the Kostka numbers as tensor product multiplicities in the BGG category $\mathcal{O}$ of $\mathfrak{s l}_{n}(\mathbb{C})$. Note that
\begin{equation}\label{h-s-expansion}
h_{\tau}=\sum_{\lambda} K_{\lambda, \tau} s_{\lambda},
\end{equation} 
which is equivalent to the definition \eqref{def} of the Schur polynomial,
see \cite{fulton1997young}*{Section 6.1}.
On the other hand, since the Schur polynomial $s_{\lambda}$ is the character of the finite-dimensional irreducible polynomial representation $V_{\lambda}$ (the Schur module) of $GL_{n}(\mathbb{C})$, by restricting to the special linear group $SL_{n}(\mathbb{C})$ and differentiating, it is also the formal character of the simple module $L(\lambda)$ of the special linear Lie algebra $sl_{n}(\mathbb{C})$, namely, \begin{equation} 
\operatorname{ch} L(\lambda)=s_{\lambda}, \text{ for } \lambda \text{ dominant and integral, i.e., } \lambda \text{ is a partition.}
\end{equation}
It follows that the Kostka number $K_{\lambda, \tau}$ is also the weight multiplicity of $\tau$ in $L(\lambda)$, namely, 
\begin{equation} 
K_{\lambda, \tau}=\operatorname{dim}\left(L(\lambda)_\tau\right),
\end{equation}
where $\lambda$ is a partition, regarded as a dominant weight of $\mathfrak{s l}_{n}(\mathbb{C})$, and $\tau$ is any sequence $\left(\tau_{1}, \ldots, \tau_{n}\right) \in \mathbb{Z}^{n}$ with $|\lambda|=|\tau|$.
By the Weyl character formula \eqref{WK} and the characters of the Verma modules \eqref{chac-verma}, we have 
\begin{equation} \label{Kostka}
\begin{aligned}
K_{\lambda, \tau}=&\operatorname{dim}\left(L(\lambda)_\tau\right)\\
=&\sum_{v \in S_{n}}(-1)^{\ell(v)}\operatorname{dim}\left(\Delta((v \cdot \lambda)_\tau\right)\\
=&\sum_{v \in S_{n}}(-1)^{\ell(v)} \mathfrak{p}(v \cdot \lambda-\tau),
\end{aligned}
\end{equation}
where $\mathfrak{p}$ is the Kostant partition function. Recall that  $\mathfrak{p}$ is defined by 
\begin{equation} \label{Kos}
\frac{1}{\prod\left(1-e^{-\alpha}\right)}=\sum_{\mu} \mathfrak{p}(\mu) e^{-\mu}, \end{equation} 
where $\alpha$ ranges over all positive roots of $sl_{n}(\mathbb{C})$. 

To prove Theorem \ref{thm1}, we use \eqref{Kostka} to give an alternative interpretation of the Kostka numbers. 
In view of \eqref{truncated} and \eqref{h-s-expansion}, we directly consider the Kostka numbers of the form $K_{\lambda, \mu+\delta-w(\nu+\delta)}$. The following result plays a key role for our proof of Theorem \ref{thm1}. 

\begin{prop} \label{KM}
The Kostka number $K_{\lambda, \mu+\delta-w(\nu+\delta)}$ is equal to $[L(\lambda) \otimes \Delta(w\cdot\nu): L(\mu)]$, the multiplicity of irreducible composition factor $L(\mu)$ in $L(\lambda) \otimes \Delta(w\cdot\nu)$, where $\lambda$, $\mu$ and $\nu$ are partitions, regarded as dominant weights for $\mathfrak{s l}_{n}(\mathbb{C})$, and $\delta=(n-1, n-2, \ldots, 1,0)$. In particular, $K_{\lambda, \mu}=[L(\lambda) \otimes \Delta(0): L(\mu)]$.
\end{prop}

\begin{proof}
One can use the \emph{standard filtration} (also known as the \emph{Verma flag}) to prove this proposition; see (\cite{humphreys2008representations}*{Theorem 3.6}). Note that $\mu$ corresponds to a dominant weight, so the multiplicity with which $\Delta(\mu)$ occurs in a standard filtration of $L(\lambda) \otimes \Delta(w\cdot\nu)$ is equal to the multiplicity of $L(\mu)$ in a Jordan-H\'older series of $L(\lambda) \otimes \Delta(w\cdot\nu)$. 

Here we give a direct computational proof using formal characters. 
If $f$ is a formal character, let $[\operatorname{ch} L(\lambda)] f$ denote the coefficient of $\operatorname{ch} L(\lambda)$ in the expansion of $f$ in terms of irreducible characters, $[\operatorname{ch} \Delta(\lambda)] f$ denote the coefficient of $\operatorname{ch} \Delta(\lambda)$ in its expansion by Verma module characters, and let $[x^{\lambda}] f$ denote the coefficient of $x^{\lambda}$ in $f$. Note that
$$[L(\lambda) \otimes \Delta(w\cdot\nu): L(\mu)]=[ \operatorname{ch} L(\mu)] \operatorname{ch} L(\lambda) \operatorname{ch} \Delta(w \cdot \nu).$$ 
Since $\mu$ is dominant, firstly we have 
\begin{align*}
[ \operatorname{ch} L(\mu)] \operatorname{ch} L(\lambda) \operatorname{ch} \Delta(w \cdot \nu) 
=&[ \operatorname{ch} \Delta(\mu)] \operatorname{ch} L(\lambda) \operatorname{ch} \Delta(w \cdot \nu)
\end{align*}
by the uni-triangularity (\ref{uni}).
Then, by the Weyl character formula (\ref{WK}), we obtain  
\begin{align*}
[ \operatorname{ch} \Delta(\mu)] \operatorname{ch} L(\lambda) \operatorname{ch} \Delta(w \cdot \nu)=&[ \operatorname{ch} \Delta(\mu)]\left(\sum_{v \in S_{n}}(-1)^{l(v)} \operatorname{ch} \Delta(v \cdot \lambda) \operatorname{ch} \Delta(w \cdot \nu)\right).
\end{align*}
Furthermore, by \eqref{chac-verma} and \eqref{Kos}, we have
\begin{align*}
[ \operatorname{ch} \Delta(\mu)]\operatorname{ch} \Delta(v \cdot \lambda) \operatorname{ch} \Delta(w \cdot \nu) =&[x^{\mu}]x^{v \cdot \lambda+w \cdot \nu} / \prod\left(1-x^{-\alpha}\right)\\
=&\mathfrak{p}(v \cdot \lambda+w \cdot \nu-\mu).
\end{align*}
Combining the above identities and (\ref{Kostka}) leads to
\begin{align*}
[ \operatorname{ch} L(\mu)] \operatorname{ch} L(\lambda) \operatorname{ch} \Delta(w \cdot \nu) =& \sum_{v \in S_{n}}(-1)^{l(v)} \mathfrak{p}(v \cdot \lambda+w \cdot \nu-\mu)=K_{\lambda, \mu+\delta-w(\nu+\delta)},
\end{align*}
as desired. In particular, taking $w=\operatorname{id}$ and $\nu=\emptyset$, we get $$K_{\lambda, \mu}=[L(\lambda) \otimes \Delta(0): L(\mu)].$$
The proof is complete.
\end{proof}

The above proposition interprets Kostka numbers as tensor product multiplicities in the BGG category $\mathcal{O}$. With this interpretation, now we are able to give a new proof of the Jacobi--Trudi identity. 

\vspace{1em}
\noindent \textit{A new proof of \eqref{GJK}.} 
On the one hand, in view of $\operatorname{ch} L(\lambda)=s_{\lambda}$ and \eqref{skew}, we have 
$$\left\langle s_{\mu / \nu}, s_{\lambda}\right\rangle=\left\langle s_{\mu}, s_{\lambda} s_{\nu}\right\rangle=[ \operatorname{ch} L(\mu)] \operatorname{ch} L(\lambda) \operatorname{ch} L(\nu)=[L(\lambda) \otimes L(\nu): L(\mu)],$$
and thus
$$s_{\mu/\nu}=\sum_{\lambda} [L(\lambda) \otimes L(\nu): L(\mu)] s_{\lambda}.$$
On the other hand, by \eqref{h-s-expansion} and proposition \ref{KM}, we have
\begin{equation*}
\begin{aligned}
\operatorname{det}\left(h_{\mu_{i}-\nu_{j}+j-i}\right)_{i, j=1}^{n}=&\sum_{w \in S_{n}} (-1)^{\ell(w)} h_{\mu+\delta-w(\nu+\delta)}\\
=&\sum_{w \in S_{n}}\sum_{\lambda}(-1)^{\ell(w)} K_{\lambda, \mu+\delta-w(\nu+\delta)}  s_{\lambda}\\
=&\sum_{w \in S_{n}} \sum_{\lambda}(-1)^{\ell(w)}
[L(\lambda) \otimes \Delta(w\cdot\nu): L(\mu)] s_{\lambda}\\
=&\sum_{\lambda} [L(\lambda) \otimes L(\nu): L(\mu)] s_{\lambda},
\end{aligned}
\end{equation*}
where the last equality follows from the Weyl character formula (\ref{WK}). This completes the proof of the Jacobi--Trudi identity \eqref{GJK}. \qed

\begin{rmk}
In particular, for the case $\nu=\emptyset$ of the Jacobi--Trudi identity, we have  \begin{equation*}
\begin{aligned}
\operatorname{det}\left(h_{\mu_{i}+j-i}\right)=&\sum_{w \in S_{n}} \sum_{\lambda}(-1)^{\ell(w)}
[L(\lambda) \otimes \Delta(w\cdot 0): L(\mu)] s_{\lambda}\\
=&\sum_{\lambda} [L(\lambda) \otimes L(0): L(\mu)] s_{\lambda}\\
=&\sum_{\lambda} [L(\lambda): L(\mu)] s_{\lambda}\\
=&s_{\mu},
\end{aligned}
\end{equation*}
which follows from the Weyl character formula (\ref{WK}) for the trivial module $L(0)$.
\end{rmk}

\begin{rmk} \label{rmk-Jacobi}
Similar arguments give a new proof of the Jacobi--Trudi-type expansion \eqref{eq-JT type} of the product of Schur polynomials: 
\begin{equation*}
\begin{aligned}
\operatorname{det}\left(h_{\mu_{i}+\nu_{n+1-j}-i+j}\right)_{i, j=1}^{n}=&\sum_{w \in S_{n}} (-1)^{\ell(w)} h_{\mu+\delta+w w_0\nu-w\delta}\\
=&\sum_{w \in S_{n}}\sum_{\lambda}(-1)^{\ell(w)} K_{\lambda, \mu+\delta-w(-w_0\nu+\delta)}  s_{\lambda}\\
=&\sum_{w \in S_{n}} \sum_{\lambda}(-1)^{\ell(w)}
[L(\lambda) \otimes \Delta\left(w\cdot(-w_0\nu)\right): L(\mu)] s_{\lambda}\\
=&\sum_{\lambda} [L(\lambda) \otimes L(-w_0\nu): L(\mu)] s_{\lambda}\\
=&\sum_{\lambda} [L(\lambda) \otimes L(\nu)^*: L(\mu)] s_{\lambda}\\
=&\sum_{\lambda} [L(\mu) \otimes L(\nu): L(\lambda)] s_{\lambda}\\
=&s_{\mu}s_{\nu},
\end{aligned}
\end{equation*}
where the last but one equality follows from tensor-hom adjunction of finite-dimensional modules of $\mathfrak{s l}_{n}(\mathbb{C})$.
\end{rmk}

We proceed to prove Theorem \ref{thm1}. 

\vspace{1em}
\noindent\textit{Proof of Theorem \ref{thm1}.} 
By Theorem \ref{thm4}, there exists an exact sequence of $\mathfrak{s l}_{n}(\mathbb{C})$-modules 
\begin{equation}\label{BGG0}
0 \rightarrow C_{\binom{n}{2}}^{\nu} \stackrel{d_{\binom{n}{2}}}{\longrightarrow} C_{\binom{n}{2}-1}^{\nu} \rightarrow \cdots \rightarrow C_{1}^{\nu} \stackrel{d_1}{\longrightarrow} C_{0}^{\nu} \stackrel{d_0}{\longrightarrow} L(\nu) \rightarrow 0,
\end{equation}
where $C_{k}^{\nu}=\bigoplus_{w \in S_n, \ell(w)=k} \Delta(w\cdot\nu)$.
For $0 \leq k \leq \binom{n}{2}$, let $\operatorname{V}\left(\nu, k\right)$ be the image of the boundary map $d_{k}$ in the BGG resolution of $L(\nu)$, and consider the exact sequence by truncating (\ref{BGG0}):
\begin{equation} \label{trun}
0 \rightarrow C_{\binom{n}{2}}^{\nu} \rightarrow C_{\binom{n}{2}-1}^{\nu} \rightarrow \cdots \rightarrow C_{k}^{\nu} \stackrel{d_{k}}{\longrightarrow} \operatorname{V}\left(\nu, k\right) \rightarrow 0.
\end{equation}
Taking the formal characters of the exact sequence (\ref{trun}), we obtain
\begin{equation} \label{posi}
\operatorname{ch}\left(\operatorname{V}\left(\nu, k\right)\right)=\operatorname{ch}\left(\operatorname{Im}\left(d_{k}\right)\right)=(-1)^{k} \sum_{w \in S_{n}, \ell(w) \geq k}(-1)^{\ell(w)}  \operatorname{ch}\left(\Delta(w \cdot \nu)\right).
\end{equation}
Now by Proposition \ref{KM}, we have
\begin{equation*}
\begin{aligned}
g_{\mu / \nu}^{k}=&(-1)^{k} \sum_{w \in S_{n}, \ell(w) \geq k} (-1)^{\ell(w)} h_{\mu+\delta-w(\nu+\delta)}\\
=&(-1)^{k} \sum_{w \in S_{n}, \ell(w) \geq k}\sum_{\lambda}(-1)^{\ell(w)} K_{\lambda, \mu+\delta-w(\nu+\delta)}  s_{\lambda}\\
=&(-1)^{k} \sum_{w \in S_{n}, \ell(w) \geq k} \sum_{\lambda}(-1)^{\ell(w)}
[L(\lambda) \otimes \Delta(w\cdot\nu): L(\mu)] s_{\lambda}\\
=&(-1)^{k} \sum_{w \in S_{n}, \ell(w) \geq k} \sum_{\lambda}(-1)^{\ell(w)}[ \operatorname{ch} L(\mu)] \operatorname{ch} L(\lambda) \operatorname{ch} \Delta(w \cdot \nu) s_{\lambda}\\
=&\sum_{\lambda} [ \operatorname{ch} L(\mu)] \operatorname{ch} L(\lambda) \operatorname{ch} \operatorname{V}\left(\nu, k\right) s_{\lambda}\\
=&\sum_{\lambda}
[L(\lambda) \otimes \operatorname{V}\left(\nu, k\right): L(\mu)] s_{\lambda}.
\end{aligned}
\end{equation*}
Almost the same arguments give the proof of \eqref{eq-tmunu}, see also Remark \ref{rmk-Jacobi}.
The proof is complete.\qed 

\begin{rmk}
The usual representation theory approach to establish Schur positivity of a given symmetric polynomial is to show that it is the characters of some finite-dimensional polynomial representation  of the general linear groups or the Frobenius images of some representation of the symmetric groups. As it can be seen from the proof, our approach is to transform the Schur positivity problem into the positivity problem of a virtual character in the BGG category $\mathcal{O}$ of $sl_{n}(\mathbb{C})$, see \eqref{posi}. It would be interesting to know whether there are other Schur positivity problems can be handled using BGG category $\mathcal{O}$ techniques.
\end{rmk}

\textbf{Acknowledgments}

The first author would like to thank Ming Fang, Hongsheng Hu, and Nanhua Xi for useful comments. The second author is supported in part by the Fundamental Research Funds for the Central Universities and the National Natural Science Foundation of China (Nos.11522110, 11971249).

\bibliography{template}

\end{document}